\documentclass[twoside,leqno]{article}

\usepackage[letterpaper]{geometry}

\usepackage{siamproceedings}

\usepackage[T1]{fontenc}
\usepackage{amsfonts}
\usepackage{graphicx}
\usepackage{epstopdf}
\usepackage{enumitem}
\usepackage{algorithmic}
\ifpdf
  \DeclareGraphicsExtensions{.eps,.pdf,.png,.jpg}
\else
  \DeclareGraphicsExtensions{.eps}
\fi

\newsiamremark{remark}{Remark}
\newsiamremark{hypothesis}{Hypothesis}
\crefname{hypothesis}{Hypothesis}{Hypotheses}
\newsiamthm{claim}{Claim}

\usepackage{amsopn}

\begin{document}

\newcommand\relatedversion{}

\title{\Large 
On the work of Jacob Tsimerman
\relatedversion}
    \author{Jonathan Pila
    \thanks{University of Oxford (\email{jonathan.pila@maths.ox.ac.uk})}}
    %email{rginder@siam.org}, \url{http://www.siam.org}).}
    %\and Rachel Ginder\footnotemark[2]    
    %\and Mitch C. Donovan\thanks{Department of Applied Mathematics, Fictional University, Austin, TX and NASA, Johnson Space Center, Houston, TX
  %(\email{donovan@fictional.edu}, \email{mcd@nasa.gov}, \url{https://www.nasa.gov/johnson/}).}}

\date{}

\maketitle

% Copyright Statement
% When submitting your final paper to a SIAM proceedings, it is requested that you include
% the appropriate copyright in the footer of the paper.  The copyright added should be
% consistent with the copyright selected on the copyright form submitted with the paper.
% Please note that "20XX" should be changed to the year of the meeting.

% Default Copyright Statement

%%%%%%%%%%%%%%%%%%%%%%%

%\fancyfoot[R]{\scriptsize{Copyright \textcopyright\ 20XX by %SIAM\\
%Unauthorized reproduction of this article is prohibited}}

%%%%%%%%%%%%%%%%%%%%%%%%%%

% Depending on which copyright you agree to when you sign the copyright form, the copyright
% can be changed to one of the following after commenting out the default copyright statement
% above.

\fancyfoot[R]{\scriptsize{Copyright \textcopyright\ 2026\\
Copyright for this paper is retained by the author}}

%\fancyfoot[R]{\scriptsize{Copyright \textcopyright\ 20XX\\
%Copyright retained by principal author's organization}}

%\pagenumbering{arabic}
%\setcounter{page}{1}%Leave this line commented out.

\begin{abstract}
\ Jacob\ \ Tsimerman\ \  was\ \   awarded  a \ 
Fields  Medal  at  the\  ICM  in
Philadelphia  on July\  23,\  2026\  \ 
``For his contribution in the recasting 
of o-minimality as a fundamental method 
of arithmetic and complex algebraic 
geometry, and his role in the proof of 
many central conjectures including
Griffiths' conjecture on the algebraicity 
of images of the period maps, and the 
André-Oort conjecture for Siegel modular 
varieties.''
We describe this work.
\end{abstract}

\section{Introduction.}

Jacob is honoured for his contributions to o-minimality 
and to resolving central problems including
conjectures of 
Griffiths and Andr\'e-Oort as described above.

The Andr\'e-Oort conjecture is a diophantine finiteness 
conjecture in the spirit of (and made in analogy with) 
the Manin-Mumford conjecture (theorem of
Raynaud \cite{RaynaudCourbes, RaynaudGeneral}) 
on torsion points on 
subvarieties of abelian varieties. It is the compositum 
of conjectures made by Andr\'e \cite{Andre89} (for 
curves in Shimura varieties) and Oort 
\cite{OortCanonicalLifts} (for subvarieties of
$\mathcal{A}_g$). Shimura varieties are a special class 
of arithmetic varieties. They have a collection of
{\it special points\/}, generalizing the CM points 
in the moduli of elliptic curves, that support a rich
arithmetic.
The Andr\'e-Oort conjecture concerns the arithmetic 
of special points on Shimura varieties.

Griffiths' algebraicity conjecture \cite{Griffiths70}, 
made in 1970, concerns a basic question arising 
from Griffiths' foundational work on variations of Hodge 
structures.
It asserts the algebraicity of the image of the 
period map, and affiliated properties,
and is a problem at the interface
of algebraic and complex geometry.

O-minimality is a branch of model theory concerned 
with tame  structures in real geometry.  It originated
in the early 80's, with van den Dries 
\cite{DriesRemarks}
exploring a decidability question of Tarski concerning 
the model theory of the real exponential function. 
The finiteness property isolated by van den Dries 
was developed in a general setting by Pillay 
and Steinhorn \cite{PillaySteinhornDefinable}. 
The crucial structure
theorem for definable sets in such a structure
was proved by Knight-Pillay-Steinhorn
\cite{KnightPillaySteinhorn}.
Geometric objects (sets and functions in/on 
$\mathbb{R}^n$) that are 
``definable in an o-minimal 
structure'' enjoy very strong tameness properties.
But showing that specific sets and functions generate
an o-minimal structure is a rather exotic and generally 
very difficult problem in real analysis; the o-minimality
of the structure we will work in here follows from work
of Gabrielov \cite{Gabrielov}, Wilkie 
\cite{WilkieModelCompleteness}, 
and van den Dries-Miller \cite{DriesMillerRanexp}.
Showing that functions of interest are definable in this 
structure is highly non-trivial as well.

The above topics are seemingly
quite distant from  each other. The Counting Theorem
of Pila-Wilkie \cite{PilaWilkie},
on the one hand, and powerful results 
on complex geometry in an
o-minimal setting, due to
Peterzil-Starchenko \cite{PeterzilStarchenkoTame}, 
on the other, provide flexible tools that
are applicable to the complex geometric objects arising
in the study of Shimura varieties and Hodge theory.
In Jacob's hands, these have proved decisive in 
obtaining results of both 
an arithmetic and geometric nature.

Much of the story of Jacob's cited work concerns the interplay
of these areas. But what I want to emphasize here
at the outset is that Jacob's contributions span 
several disparate aspects of that interaction.
Indeed, Jacob and his collaborators have made fundamental
advances in (at least) four separate aspects 
involved in applying o-minimality to 
arithmetic and Hodge-theoretic problems.
(I will not set out in detail here which work was 
done in collaboration and with whom; in the sequel 
full attributions will be given.) These are:

\smallskip

(1) proving crucial arithmetic properties of 
special points that, combined with o-minimal 
techniques and functional transcendence results,
enabled the proof of the Andre-Oort conjecture
for Siegel modular varieties, and indeed in full;

\smallskip

(2) establishing suitable functional transcendence 
results (partly via o-minimality, partly otherwise)
needed for Andr\'e-Oort but also going beyond those to 
secure a functional version of the Andr\'e-Grothendieck 
period conjecture;

\smallskip

(3) deriving definability results that bring further
objects into the o-minimal realm, crucially including the natural 
objects of variational Hodge theory;

\smallskip

(4) developing an ``o-minimal GAGA'', a counterpart to the 
classical GAGA of Serre \cite{SerreGAGA}, 
which (in conjunction with 
(3)) may be applied in the setting of Hodge theory,
and leads to a proof of the Griffiths conjecture.

\smallskip

In the sequel I will introduce each of these topics
and describe the beautiful and spectacular results
that Jacob's work has given us.
For a different perspective see
\cite{Simpson26}.

\section{Context of the Andr\'e-Oort conjecture.}

To motivate the conjecture, we first consider the wider 
diophantine context. Diophantine geometry studies, 
in the first instance, the solution of systems of 
algebraic equations in integers, or in rational numbers.
In general, the existence of integer solutions is known to
be algorithmically undecidable (recently, also for
rational points, according to OpenAI
\cite{OpenAIH10Q}). A (possibly)
more accessible question, and a convenient one for us is:

\smallskip

{\it Which complex algebraic varieties have a Zariski-dense
set of rational points?}

\smallskip

When $V$ is an algebraic curve (i.e. a 
one-dimensional complex algebraic variety) an
answer was conjectured by Mordell \cite{Mordell}
and affirmed by Faltings \cite{FaltingsMordell}: 
only a curve of genus 0 or 1 can have a 
Zariski-dense (i.e. infinite) set of rational points. 
A curve of genus 0 or 1 may or may not have
an infinitude of rational points; it will do so however
over some finite extension of $\mathbb{Q}$. If one asks:
which curves have a Zariski-dense set of points 
over some finite
extension of $\mathbb{Q}$
(i.e. a number-field), the answer 
(by Faltings \cite{FaltingsMordell}) is:

\smallskip

{\it Among curves, only those of genus $g<2$ can have a 
Zariski-dense set of rational points over a number-field.\/}

\smallskip

For varieties of dimension 2 or greater, the above 
question  is very much open. Bombieri conjectured 
that a surface of
general type does not have a Zariski-dense set of rational
points, and the Bombieri-Lang conjecture asserts that
a variety is {\it Mordellic\/} (i.e. has only finitely many 
points rational over any fixed finitely generated
field) outside its
``special set'', which is geometrically defined.
(Thus framed the question seems somewhat resistant 
to undecidability phenomena.)
Faltings \cite{FaltingsBig} 
%(see also \cite{FaltingsGeneral, OortTheGeneral}) 
answered the question above for subvarieties of 
abelian varieties:

\smallskip

{\it Among subvarieties of abelian varieties, 
only abelian subvarieties and their translates
can have this property.}

\smallskip

An abelian variety $A$ defined over $\mathbb{C}$ may be 
realized as the quotient of $\mathbb{C}^n$, 
where $n=\dim A$, by a lattice. This makes the 
abelian group structure on $A$ (inherited from 
addition in $\mathbb{C}^n$) transparent. 
Division points of the lattice give torsion
points in $A$, which are thus seen to be 
analytically dense. The Manin-Mumford conjecture addresses
the following question:

\smallskip

{\it Let $A$ be an abelian variety.
Which subvarieties $V\subset A$ have a 
Zariski-dense set of torsion points?}

\smallskip

The answer is due to Raynaud \cite{RaynaudCourbes, 
RaynaudGeneral}:

\smallskip

{\it Such $V$ are precisely the ``torsion cosets''
i.e. the translates of abelian subvarieties 
by torsion points.}

\smallskip

Note that our passage to Manin-Mumford via 
Faltings ``Big'' theorem, while convenient for
our presentation, 
is anhistorical. The Manin-Mumford conjecture was 
originally formulated in the 60's for curves 
embedded in their Jacobians. Lang subsequently formulated the 
``Mordell-Lang'' conjecture, unifying and generalizing 
the Mordell and Manin-Mumford conjectures, 
which was proved later by Faltings (in conjunction 
with results of, and with more general versions proved by,
Hindry, Laurent, McQuillen, Raynaud, and Vojta)
see \cite{OortTheGeneral}.

\section{The Andr\'e-Oort conjecture.}

{\it Shimura varieties\/} are certain very special
quasi-projective varieties, (canonically) defined over
$\overline{\mathbb{Q}}$,
that arise as moduli spaces (of suitable Hodge structures)
and are central objects in arithmetic geometry and the
Langlands programme. Initial examples were defined and
studied by Shimura 
\cite{ShimuraAnalyticFamilies, ShimuraConstruction} and
the modern definition is due to Deligne \cite{Deligne71}.
A convenient reference for a precise definition is
\cite{Daw15}.

The paradigm examples are the {\it Siegel modular varieties\/}
$\mathcal{A}_g, g=1,2, \ldots$, the moduli space of
principally polarized complex abelian varieties of 
dimension $g$, studied by Siegel \cite{SiegelModularVariety}. 
Already $\mathcal{A}_1=\mathbb{A}^1$,
arising as the image of the modular function 
$j: \mathbb{H}\rightarrow \mathbb{C}$
(in which guise it is also known as the 
{\it modular curve\/} $Y(1)$), which parameterizes
elliptic curves defined over $\mathbb{C}$ by their
$j$-invariant, has a rich and beautiful arithmetic theory
revolving around its ``special points''.

A Shimura variety is endowed with a collection
of {\it special points\/}. These $\overline{\mathbb{Q}}$ points
are the moduli corresponding to the objects with
``most symmetries''.
The special points of $\mathcal{A}_1$ are 
the $j$-invariants of elliptic curves possessing
endomorphisms other than
multiplication by $m\in \mathbb{Z}$; such elliptic curves
are said to have
{\it complex multiplication\/}, and their $j$-invariants
are classically known as {\it singular moduli\/},
and intimately related with class-field theory.
An elliptic curve $E$ defined over $\mathbb{C}$ 
may be realized as 
the quotient $\Lambda\backslash\mathbb{C}$
of $\mathbb{C}$ by a lattice, which may be
normalized to be of the form $\Lambda=\mathbb{Z} + 
\mathbb{Z}\tau$, where $\tau\in\mathbb{H}$.
The $j$-invariant of $E$ is then $j(E)=j(\tau)$ and is
a singular modulus precisely when 
$[\mathbb{Q}(\tau):\mathbb{Q}]=2$.
Quadratic points are analytically dense in 
$\mathbb{H}$, and consequently singular moduli 
are analytically dense in $\mathcal{A}_1$.

In general, a Shimura variety $S$ may also be realized 
as the quotient of a suitable Hermitian symmetric space 
$\mathbb{H}_S$ by an arithmetic subgroup of its group of
biholomorphisms, and has an analytically dense set of
special points.
In the case $S=\mathcal{A}_g$, the parameterized objects
are abelian varieties of dimension $g$, and the
``special'' ones are the CM (complex multiplication)
abelian varieties. For simple $A$, this means that the
endomorphism ring of $A$ is an order in a CM field
of degree $2g$ over $\mathbb{Q}$
(and so is as large as is possible); 
in the non-simple case some care must be taken.
The Andr\'e-Oort conjecture addresses the question:

\smallskip

{\it Let $S$ be a Shimura variety. 
Which subvarieties $V\subset S$ have a 
Zariski-dense set of special points?}

\smallskip

In addition to special points, a Shimura variety $S$ has 
a broader collection of {\it special subvarieties\/}. 
These are ``Shimura subvarieties'' whose
structure as a Shimura variety is compatible with that
of $S$ (see \cite{Daw15}).
The special points are precisely the special
subvarieties of dimension zero.
If positive dimensional, a special subvariety
$T\subset S$ has a Zariski-dense (and even analytically dense)
set of special points of $S$ on it. The Andr\'e-Oort
conjecture predicts, in analogy with Manin-Mumford, that
these are the only such subvarieties.
Let us mention here that varieties will be identified
with their sets of complex points, and subvarieties 
are assumed closed in their ambient.

\medskip
\noindent
{\bf 3.1. Conjecture.\/} 
(Andr\'e-Oort \cite{Andre89, OortCanonicalLifts}) 
Let $S$ be a Shimura
variety. A subvariety $V\subset S$ has a Zariski-dense 
set of special points iff it is a special subvariety.

\medskip

The Andr\'e-Oort conjecture is trivial for one-dimensional
Shimura varieties. The first non-trivial result towards 
it was proved by Andr\'e \cite{Andre98}, for a curve in
$S=\mathcal{A}_1^2=(\mathbb{A}^1)^2$. 

\medskip
\noindent
{\bf 3.2. Theorem.\/} (Andr\'e \cite{Andre98})
{\it The Andr\'e-Oort conjecture holds
for $\mathcal{A}_1^2$.}

\medskip

The special points in $\mathcal{A}_1^2$
are pairs $(\sigma_1, \sigma_2)$ of singular moduli.
In addition to the special points, and $S$, the special 
subvarieties are: ``vertical'' lines of the form
$\{(\sigma, y): y\in \mathbb{A}^1\}$ where $\sigma$
is a singular modulus, the ``horizonal'' lines 
$\{(x, \sigma): x\in \mathbb{A}^1\}$, 
and the modular curves $\Phi_N(x,y)=0$
defined by the classical modular polynomials. 
These are the images under the map 
$\mathbb{H}^2\rightarrow 
\mathbb{C}^2, (u,v)\mapsto (j(u), j(v))$
of curves of the form $v=gu$ where 
$g\in{\rm GL}_2(\mathbb{Q})^+$ acts by fractional 
linear transformations
(and ``+'' indicates positive determinant).

Around the same time, Edixhoven \cite{EdixhovenTwo} 
gave a  different proof of the same result, 
conditional on GRH,
and later generalized it to $\mathcal{A}_1^n$
\cite{EdixhovenProducts}.
This line of attack was further pursued in a series 
of papers leading
a full affirmation of the conjecture, conditional on GRH,
by Klingler, Ullmo, and Yafaev 
\cite{KlinglerYafaevAO, UllmoYafaevAO}.
In these developments, equidistribution results for
special subvarieties played a crucial role
\cite{ClozelOhUllmo, ClozelUllmo, CohenHyperbolic,
ZhangAO}. See \cite{YafaevSurvey}.

A different approach, using the theory of ``o-minimal 
structures'', enabled unconditional results and 
led to Jacob's proof of 
the Andr\'e-Oort conjecture for $\mathcal{A}_g$.

\section{O-minimality.}

The notion of o-minimality emerged from the model-theoretic 
study of the real exponential function. 

Recall that a {\it structure\/}, in the sense of 
model theory, consists of a set $M$, the ``universe''
of the structure, together with some specified 
``structure'' on $M$ comprising: some constants $c_i\in M$, 
for $i\in I$,
functions $f_j: M^{n(j)}\rightarrow M, j\in J$ 
(of arities $n: J\rightarrow \mathbb{N}_{\ge 1}$), 
and relations $R_k\subset M^{m(k)}, k\in K$
(of arities $m: K\rightarrow\mathbb{N}_{\ge 1}$).
(The index sets $I, J, K$ are arbitrary, possibly empty.)
One writes
$$
\mathcal{M}=(M, c_i, f_j, R_k\}
$$
to denote this structure (where I have suppressed the 
index sets).
For example, $(\mathbb{R}, <, +, \times, 0, 1)$ denotes
the structure of the real numbers as an
ordered field. 

Corresponding to a structure is a (first-order)
language comprising symbols for constants, functions, 
and relations
with the same index sets and arities. Formally, 
one distinguishes the symbols of the language from their
interpretation in a structure for the language,
but this is usually honoured in the breach.
In addition to symbols corresponding to the 
structure, the language has variable symbols
$x,y,\ldots, x_1, x_2,\ldots$, logical symbols 
$\wedge, \lor, \lnot,
\forall, \exists$, and $=$, and parentheses (for
unique readability).

The formulas of the language are defined syntactically,
but can be informally described as the ``meaningful 
strings'' in the language. Variables appearing may or 
may not be bound by a quantifier ($\forall$ or $\exists$),
formulas with no free variables are called 
{\it sentences\/}. A sentence $\phi$
then can be read as an assertion that may or 
may not hold in a given structure $\mathcal{M}$
for the language. If it does, one writes 
$\mathcal{M}\models \phi$.

If $\phi(x_1,\ldots, x_n)$ is a formula in which
the free variables are $x_1,\ldots, x_n$, then the
set 
$$
Z=\{x=(x_1,\ldots, x_n)\in M: \mathcal{M}\models \phi(x)\}
$$
is a {\it definable set\/} in the structure $\mathcal{M}$.
On the basic setup of model theory see \cite{MarkerGTM}.

Tarski \cite{Tarski} proved that the structure
$(\mathbb{R}, <,  +, \times, 0, 1)$ is decidable;
i.e., he gave a decision procedure for the 
sentences. 
Note that the inclusion of $<$ is unnecessary here, 
as the relation $x < y$ is definable in 
$(\mathbb{R}, +, \times, 0, 1)$ by the formula
$\exists t(\lnot t= 0 \land y=x+t^2)$, so its
inclusion does not change what can be expressed
(i.e. same definable sets).
But o-minimal structures are {\it by definition\/} 
expansions of ordered
structures, so we prefer to explicitly include it.

Tarski's decision procedure determines
whether a given sentence holds in $\mathbb{R}$ where
the symbols $+, \times, 0, 1$ are interpreted in 
the usual way. The complex field 
$(\mathbb{C}, +, \times, 0, 1)$ is also decidable, 
but this is easier: by Chevalley's theorem such
decisions resolve to questions of ideal membership. The
rational number field is however undecidable 
\cite{Robinson}. This is an easier problem than 
Hilbert's 10th problem for the rational numbers 
as sentences, which can involve 
quantifiers, are much richer than systems of 
diophantine equations.

Tarski \cite{Tarski} asked whether the structure 
$\mathbb{R}_{\exp}=(\mathbb{R}, <, +, \times, 0, 1, \exp)$
is decidable, where $\exp$ denotes the real
exponential map $\mathbb{R}\rightarrow\mathbb{R}$.
By contrast, the structure 
$\mathbb{R}_{\sin}=(\mathbb{R}, <, +, \times, 0, 1, \sin)$
is undecidable, as it defines the integers
(as does the structure $(\mathbb{C}, +, \times, 0, 1, \exp)$
where $\exp$ is the complex exponential).

Dissecting the various strands in Tarski's work 
on $(\mathbb{R}, <, +, \times, 0, 1)$
in the context of this question, van den Dries
\cite{DriesRemarks} notes that one aspect that might
generalize is the tame nature of definable sets,
specifically in the form of admitting ``cylindrical 
decomposition'' as introduced by the work of Collins 
\cite{Collins}
on decision algorithms for the real field. 
van den Dries frames the crucial condition which
becomes the key defining property of an o-minimal
structure, calling it ``finite type'',
and expects this property to hold for
(the theory of) $\mathbb{R}_{\exp}$.

The next step was taken by Pillay and Steinhorn 
\cite{PillaySteinhornDefinable, 
PillaySteinhornDefinable1} who
explored the ``finite type'' condition in a general 
setting, and named the structures ``o-minimal''
(i.e. ``order-minimal'': they are the analogues of
``minimal'' structures in the presence of an order).
Together with Knight \cite{KnightPillaySteinhorn}
they proved the crucial Cell Decomposition Theorem.
This proves that the notions ``o-minimal'' and ``strongly 
o-minimal'' coincide. (By contrast, minimal
structures do not coincide with strongly minimal ones;
the latter have played a key role in ``stable''
model-theory, e.g. in Morley's proof of the Los
conjecture.) The very strong tameness properties 
of definable sets in o-minimal structures flow 
from this theorem (e.g. definable functions are 
continuous, and even differentiable, 
except for finitely many points, uniformly in families,
etc. etc. see e.g. \cite{DriesTT, DriesMillerCategories}).

The setting of o-minimality is a structure $\mathcal{R}$
that  {\it expands\/} a dense linear order 
without endpoints $(R, <)$. 
This means that the structure on $\mathcal{R}$ includes
an order, under which its universe $R$ is densely linearly
ordered, with no largest or smallest elements, as well as
perhaps further structure. Such $\mathcal{R}$ is 
usually indicated as $\mathcal{R}=(R, < ,\ldots)$.

\medbreak
\noindent
{\bf 4.1. Definition.\/}
Let $\mathcal{R}=(R, < ,\ldots)$ be an expansion of a dense
linear order $(R, <)$ without endpoints. Then $\mathcal{R}$
is {\it o-minimal\/} if every definable
subset of $R$ is a finite union of points and open 
(possibly unbounded) intervals.

\medbreak

Note that, in o-minimality, it is standard to consider
definable sets as being those definable ``with parameters'', 
meaning that all constants from the universe may be used as 
parameters in defining sets, and we do so in the above 
definition. Thus, in the structure
$(\mathbb{R}, <, +, \times, 0, 1)$ of the real field, 
the  definable sets in this sense are all 
semi-algebraic sets,
while under the strict definition of definability 
(using only constants in the structure)
one gets only semi-algebraic sets defined over
the algebraic closure of $\mathbb{Q}$ in $\mathbb{R}$.
The underlying order in applications is usually 
$(\mathbb{R}, <)$, but it is important for the theory to
take place in the more general setting.

To be interesting and useful, however, there should be
rich and flexible o-minimal structures beyond the
prototype afforded by $(\mathbb{R}, <,  +, \times, 0, 1)$
itself, while the primary target of the theory was 
the structure $\mathbb{R}_{\rm exp}$.
A first example, already very rich, is afforded by the
structure of {\it restricted analytic functions\/},
usually denoted 
$\mathbb{R}_{\rm an}=(\mathbb{R}, <, +, \times, 0, 1, f_i, 
i\in I)$ having a symbol $f_i$ for every function
$f_i:B\rightarrow\mathbb{R}$ where 
$B\subset \mathbb{R}^n$, for some $n\ge 1$, is a 
compact box 
and $f_i$ extends to a real analytic function on some open
neighbourhood of $B$. In the structure, one includes the
restriction of $f_i$ to $B$ (setting it to be 0 outside).
van den Dries
\cite{DriesTarskiSeidenberg} observes that the 
o-minimality of $\mathbb{R}_{\rm an}$
is a consequence a famous theorem of Gabrielov 
\cite{Gabrielov}, asserting that 
the complement of a subanalytic set $A\subset \mathbb{R}^n$
is subanalytic.

This structure includes the restriction of $e^x$ to any
compact interval, but (by a theorem of 
Lojasiewicz) does not include the full graph of the
real exponential.
The o-minimality of $\mathbb{R}_{\rm exp}$ was proved by
Wilkie \cite{WilkieModelCompleteness} 
(who indeed proves that it is model-complete, i.e. that every definable set is the projection of a 
quantifier-free-definable set, 
from which o-minimality follows by Khovanskii's theory 
\cite{Khovanskii} of ``Fewnomials''). 

Combining these structures, van den Dries and Miller
\cite{DriesMillerRanexp} show that the structure 
$\mathbb{R}_{\rm an\ exp}$ that expands 
$\mathbb{R}_{\rm an}$ by including the full graph 
of $e^x$ is o-minimal; a second proof is given in
\cite{DriesMM} (and later it was shown that 
the ``Pfaffian closure'' of any o-minimal structure 
over $\mathbb{R}_{\rm alg}$ is again o-minimal: 
\cite{Speissegger},
and this always includes $\exp$).

\medbreak
\noindent
{\bf 4.2. Theorem.} (van den Dries-Miller 
\cite{DriesMillerRanexp})
{\it The structure $\mathbb{R}_{\rm an\ exp}$ is o-minimal.}
\medbreak

To date, the structure $\mathbb{R}_{\rm an\ exp}$
has been adequate for all the 
applications to Shimura varieties 
and Hodge theory. For some other structures 
(for example in which functions such as $\Gamma(x), 
\zeta(x)$, suitably restricted, are definable) see 
\cite{DriesSpeissegger}.
%for an application
%of the structure of multisummable Gevrey series
%of \cite{DriesSpeissegger} to a Manin-Mumford-type 
%problem involving zeros of Bessel functions see
%\cite{BinyaminiKiroKozmaPila}.
Further, there exist incompatible o-minimal structures
(i.e. two o-minimal structures whose union generates a 
structure that fails to be
o-minimal \cite{RolinSpeisseggerWilkie}).
Thus while o-minimality is sometimes cited as a 
realization of Grothendieck's vision 
\cite{Grothendieck} of ``topologie moder\'ee'', it
does not afford one ``canonical'' such realization.
The standard reference for o-minimality is the book
\cite{DriesTT} by van den Dries. Many further useful 
properties are gathered in 
\cite{DriesMillerCategories}.

\section{Point-counting in 
o-minimal structures.}

O-minimality connects with the Andr\'e-Oort
conjecture via counting results
\cite{PilaWilkie} for rational points on definable sets 
in o-minimal structures 
and their application to diophantine problems.
%via the  Pila-Zannier strategy introduced 
%in \cite{PilaZannier}.
O-minimal structures provide a natural setting for
point-counting results on tame real sets.

Let $Z\subset \mathbb{R}^n$. We want to study
the rational points on $Z$, which will in general be
non-(semi-)algebraic. If, say, $Z$ is some transcendental
real-analytic set, we might expect rational points on it
to be rather sparse. Though examples such as the graph of
$y=2^x$ show that they may be infinite. So we might
seek results saying that points are sparse in terms of
how their number grows with their height. We further
observe that
a transcendental real-analytic set of (real) dimension 
at least 2 might contain positive-dimensional 
semi-algebraic sets, such as lines or line segments, 
that might contain rather many rational points in
a height-density sense. E.g. the graph of $z=x^y$
on $[1,2]^2$, say.
In analogy with the idea of the
``special set'' in diophantine geometry, a reasonable
target result could be to show that, on suitably tame
real sets, rational points are sparse outside
such semi-algebraic subsets.
To state a precise result, we make some definitions.

\medskip
\noindent
{\bf 5.1. Definition.\/} The {\it height\/} of a rational
number $q=a/b$, where $a,b\in \mathbb{Z}, b\ge 1$, and 
${\rm gcd}(a,b)=1$ is set to be $H(q)=\max(|a|, b)$.
The {\it height\/} of an $n$-tuple 
$x=(x_1,\ldots, x_n)\in\mathbb{Q}^n$ is defined to 
be $H(x)=\max H(x_i)$

\medskip
\noindent
{\bf 5.2. Definition.\/} The {\it algebraic part\/}
of a set $Z\subset\mathbb{R}^n$ is defined to be the
union of all connected positive-dimensional
semi-algebraic subsets $A\subset Z$.
Its complement $Z-Z^{\rm alg}$ is the 
{\it transcendental part\/} of $Z$.

\medskip
In its simplest form, the point-counting result will 
simply exclude the algebraic part of $Z$. The algebraic 
part may be considered as a coarse analogue of the
special set in diophantine geometry.

\medskip
\noindent
{\bf 5.3. Definition.\/} The {\it counting function\/}
for $Z$ is the function
$$
N(Z,H)=\#\{x\in Z: x\in \mathbb{Q}^n, H(x)\le H\}.
$$

\medskip
\noindent
{\bf 5.4. Theorem.\/} (Pila-Wilkie \cite{PilaWilkie}).
{\it Let $Z\subset\mathbb{R}^n$ be a definable set
in an o-minimal structure, and $\epsilon>0$. Then
$$
N(Z-Z^{\rm alg}, H) \le c(Z, \epsilon) H^\epsilon
$$
for some positive constant $c(Z, \epsilon)$.}

\medskip

This result can be refined in various ways. First,
the implied constant is uniform over any {\it definable
family\/} of sets. This just means sets of the form 
$Z_t, t\in\mathbb{R}^m$ that are fibres of
the projection to $\mathbb{R}^m$ of some definable
set $Z\subset \mathbb{R}^m\times\mathbb{R}^n$.
This uniformity, which flows from the Cell 
Decomposition Theorem for definable sets 
in o-minimal structure, is crucial in the proof, at the heart
of which is a parameterization theorem for definable sets
generalizing such a result for semi-algebraic sets
due to Yomdin and Gromov 
\cite{Gromov, YomdinVolume, YomdinResolution}.

Second, one has an estimate of the same form for 
algebraic points of fixed degree \cite{PilaAlgebraic}: 
here one sets
$$
N(Z, d, H)=\#\{x\in Z: [\mathbb{Q}(x_i):\mathbb{Q}]\le d, 
H(x_i)\le H, i=1,\ldots, n\}
$$
where  $H$ denotes the multiplicative Weil height
on $\overline{\mathbb{Q}}$, and
one has, for a definable set $Z$,
$$
N(Z-Z^{\rm alg}, d, H) \le c(Z, d, \epsilon)H^\epsilon.
$$

Third, while in this form the algebraic part is excluded,
one can state a more refined version that includes it,
essentially showing that, up to height $H$, the rational
points are accounted for by $c(Z, d, \epsilon)H^\epsilon$
``components'' of the algebraic part, and points
\cite{PilaAlgebraic}.

And fourth, while in a general o-minimal structure, 
the counting estimate is ineffective and
cannot be substantially improved (see \cite{PilaQJM}), 
if one restricts to more ``controlled'' o-minimal 
structures,
e.g. where the basic generating functions satisfy suitable
differential equations, one has results in which
the upper bound is effective and/or poly-logarithmic
\cite{Binyamini22, BinyaminiNovikovZak, 
BinyaminiLogNoetherian}.

The counting theorem provides an initial bridge between 
o-minimality and arithmetic problems, via the
so-called Pila-Zannier strategy, introduced
in \cite{PilaZannier} to give a new proof
of the Manin-Mumford conjecture, and subsequently
efficacious in various cases e.g.
\cite{MasserZannierTorsionAnomalous, MasserZannierPell}
of the Zilber-Pink conjecture 
\cite{Zilber, BMZ07, PinkPreprint}.
On this conjecture generally see
\cite{PilaBook, ZannierBook}.
In the next section we set out how this strategy approaches
the Andr\'e-Oort context together with the 
ingredients required for it to succeed.

\section{O-minimality and the Andr\'e-Oort conjecture.}

The uniformization $u:\mathbb{H}_S\rightarrow S$ of a 
Shimura variety $S$ by a Hermitian symmetric domain 
$\mathbb{H}_S$, invariant under a suitable arithmetic 
group $\Gamma$ of automorphisms of $\mathbb{H}_S$,
provides the basic elements to which the point-counting 
strategy are applicable: one may realize $\mathbb{H}_S$
in its Harish-Chandra embedding as an open semi-algebraic 
subset in $\mathbb{C}^n, n=\dim S$, giving it a natural
$\overline{\mathbb{Q}}$-structure, in which the 
pre-images of special points are algebraic and of 
some bounded degree (indeed, each is the unique fixed
point in $\mathbb{H}_S$ of some $g\in G(\mathbb{Q})$
where $G$ is the governing real algebraic $\mathbb{Q}$-group
(e.g. $G={\rm Sp}_{2g}$ for $S=\mathcal{A}_g$);
\cite{DawOrr16}, p5.

In the sequel it is convenient to talk about
``algebraic subvarieties $V\subset \mathbb{H}_S$'' 
although, strictly speaking, there are none of
positive dimension.
For example, in the simplest case $S=Y(1)^n$ we 
have $\mathbb{H}_S=\mathbb{H}^n$ where $\mathbb{H}$
is the complex upper half-plane. A positive-dimensional
algebraic subset of $\mathbb{C}^n$ will be dominant
to some factor, and so cannot be contained in 
$\mathbb{H}^n$. So we will specify that an {\it algebraic 
subvariety of $\mathbb{H}_S$\/} will mean a connected
component of $\mathbb{H}_S\cap V$ for some (honest)
algebraic subvariety $V\subset\mathbb{C}^n$.

We must also introduce now the notion of a 
{\it weakly special subvariety\/}. These 
include the special subvarieties, and may be
characterized, following 
\cite{UllmoYafaevCharacterization}, as the
subvarieties $W\subset S$ such that some component
(equivalently, every component) of $u^{-1}(W)$
is algebraic (in the above sense). They are the
``bi-algebraic'' subvarieties. They can also be
characterized as the totally geodesic subvarieties
\cite{Moonen}. A weakly special subvariety is special
iff it contains a special point.

Ullmo \cite{UllmoQuelques}
sets out four ingredients that are sufficient to
successfully apply the point-counting strategy 
to establish
the Andr\'e-Oort conjecture for a given Shimura variety
$S$.

\smallskip

1. Definability of $u$ on some fundamental domain
$F_S\subset\mathbb{H}_S$ for the action of $\Gamma$.

\smallskip

2. A suitable upper bound on the height of the pre-image
in $F_S$ of a special point $P\in S$.

\smallskip

3. ``Ax-Lindemann'': the Zariski closure of $u(A)$,
for an algebraic subvariety $A\subset \mathbb{H}_S$,
is weakly special.

\smallskip

4. A suitable lower bound for the size of the Galois 
orbit of a special point.

\smallskip

Items 2 and 4 are expressed in terms of a parameter
$\Delta(P)$ measuring the ``complexity'' of a 
special point. More about it later, but
Galois conjugate special points have the same 
complexity. One then
wants, in 2, that the multiplicative height $H(Q)$,
where $Q\in F_S$ with $u(Q)=P$, satisfies 
$H(Q) \le c\Delta(P)^\delta$ for some positive
constants $c, \delta$ depending only on $S$, 
while, in 4, one wants
$$
[\mathbb{Q}(P):\mathbb{Q}] \ge c'\Delta(P)^{\delta'}
$$
for some positive constant $c', \delta'$ depending only on $S$.
%(Here we have notationally distinguished the constants,
%but in general we will assume that constants are fresh
%on each occasion.)

Now suppose that $V\subset S$ is an algebraic subvariety.
We may assume that $V$ is defined over a number field
$K$ (note that Shimura varieties are (canonically)
defined over a number field).
Indeed, all points of $V$ defined over $\overline{\mathbb{Q}}$
lie in a subvariety $V'\subset V$ defined over
a number field.
We set $Z=u^{-1}(V)\cap F_S$. By 1, this is a definable set.
Suppose that $V$ contains a special point $P$.
Then $P$ has at least 
$c'\Delta(P)^{\delta'}/[K:\mathbb{Q}]$ conjugates $P'$
that also lie in $V$, and are special (with same
complexity $\Delta$). The pre-images under $u$ in $F_S$
of these are of some bounded degree $d=d(S)$ and
have height at most $c\Delta(P)^\delta$,
and there are at least $c'\Delta^\delta$ of them.
If we set $H=c\Delta^\delta$ we thus have at least
$c''H^{\delta/\delta'}$ points in $Z$ of some fixed
degree and height $H(Q)\le H$. 

Applying the Counting Theorem with some 
$\epsilon<\delta/\delta'$ we find that, once
$\Delta$ is sufficiently large, most (and we will 
see in fact all) of the conjugates $P'$ of $P$ 
over $K$ lie in some positive dimensional special 
subvariety $T\subset V$.
To see this we need to use item 3, Ax-Lindemann. 
By the Counting Theorem, most of the $Q'$ will lie in the 
algebraic part of $Z$. 
This is a union of positive dimensional
{\it real\/} semi-algebraic subsets $A\subset Z$. 
But since $u^{-1}(V)$ is complex analytic, such $A$ is 
contained (at least locally at any smooth point) in 
some complex algebraic subvariety $B\subset u^{-1}(V)$, 
in the sense set out above.

By Ax-Lindemann, the Zariski closure of $u(B)$ is weakly 
special, and it is a subvariety of $V$. So we find
that ``most'' conjugates $P'$ of $P$ over $K$ lie 
in a positive-dimensional weakly special subvariety
$T\subset V$. And as $P'$ is special we see that $T$ is
special. And since special subvarieties are preserved 
under Galois conjugation, we conclude that all
conjugates of $P$ over $K$ are contained in some
positive-dimensional special subvariety.

This already shows that $V$ contains only finitely many 
special points that are not contained in a larger special
subvariety contained in $V$ and the argument is concluded 
by induction. For the induction one shows
that the maximal weakly special
subvarieties contained in $V$ come in finitely many families,
and are parameterized by a subvariety of a suitable
parameter space which is itself a Shimura variety.
This is not a difficult point but we won't dwell on it.

Let me set out next the work establishing items 1 and 2 
in the context of Shimura varieties before expanding on 
Jacob's contributions to items 3 and 4.
For this one needs a suitable fundamental set (or domain)
for the $\Gamma$ action on $\mathbb{H}_S$ and in the case
of Siegel modular varieties the natural choice
is the Siegel fundamental set $F_g$ introduced by 
Siegel \cite{SiegelModularVariety}.

In the simplest case of $Y(1)$, and its cartesian powers,
the definability amounts to showing that the $j$-function 
is definable on the standard fundamental domain for the
action of ${\rm SL}_2(\mathbb(Z)$ on the upper half-plane.
As the fundamental domain lies within a vertical strip,
and is bounded away from the real axis, the definability
in $\mathbb{R}_{\rm an\ exp}$ is a consequence of
the $q$-expansion, which converges for $|q|\le c<1$,
as observed in \cite{PeterzilStarchenkoUniform}.

For higher Siegel modular varieties, the definability requires 
a delicate study of the $\theta$ functions giving the
projective embedding, and was carried out
by Peterzil-Starchenko \cite{PeterzilStarchenkoDefinability}. 

%The general case of arithmetic quotients 
%(including the general case of Shimura varieties) 
%was done by Klingler-Ullmo-Yafaev
%\cite{KlinglerUllmoYafaevHyperbolicAL}. This required a 
%careful analysis of the machinery of toroidal 
%compactification.
%
%For the following theorem, we consider an arithmetic 
%variety $X$, arising as a quotient 
%$S=\Gamma\backslash\mathbb{U}$ of a hermitian symmetric
%domain. We assume $\mathbb{U}\subset\mathbb{C}^n$
%in its Harish-Chandra realization. Here $\Gamma$ is an
%arithmetic subgroup of the group of
%biholomorphisms of $\mathbb{U}$.
%
%\medbreak
%\noindent
%{\bf Theorem.} (Klingler-Ullmo-Yafaev 
%\cite{KlinglerUllmoYafaevHyperbolicAL})
%{\it Let $u: \mathbb{U} \rightarrow S$ be the uniformizing
%map of an arithmetic variety. There exists a semi-algebraic
%Siegel set $\Sigma\subset \mathbb{U}$ such that the 
%restriction 
%$u\vert_\Sigma: \rightarrow S$ of the uniformizing map
%is definable in $\mathbb{R}_{\rm an\ \exp}$.}

\medbreak
\noindent
{\bf 6.1. Theorem.\/}
(Peterzil-Starchenko \cite{PeterzilStarchenkoDefinability})
{\it There is an open set $U \subset \mathbb{H}_g$,
containing $F_g$, such that the restriction 
$u: \mathbb{H}_g\rightarrow\mathcal{A}_g$
to $U$ is definable in the o-minimal structure 
$\mathbb{R}_{\rm an\ exp}$.\/}

\medbreak

The theorem is proved via a corresponding statement
for the  universal family of abelian varieties of any 
given fixed polarization type (using theta functions
in two sets of variables $\theta
\bigl[ \begin{smallmatrix}a \\ 
 b \end{smallmatrix} \bigr]
(z, \tau)$, with characteristics).
For a general Shimura variety, the corresponding definability
statement is due to  Klingler-Ullmo-Yafaev
\cite{KlinglerUllmoYafaevHyperbolicAL}. This required a 
careful analysis of the machinery of toroidal 
compactification.
For mixed Shimura varieties the result is due to 
Gao \cite{GaoTowardsMixed}.

Moving to item 2, our complexity measure $\Delta(x)$
for a special point
in $\mathcal{A}_g$ will be the absolute value of the
discriminant of the center of the endomorphism ring
of the corresponding abelian variety.
In the case of $\mathcal{A}_1=Y(1)$ a suitable bound
can be extracted elementarily from the proof (as given
e.g. in Serre \cite{SerreArithmetic}) that the 
standard fundamental domain is indeed a fundamental 
domain. Already for the Siegel modular varieties 
it is quite non-trivial, and due
to Jacob in his PhD Thesis, appearing in 
\cite{PilaTsimerman13}.

\medbreak
\noindent
{\bf 6.2. Theorem.} (Tsimerman \cite{PilaTsimerman13})
{\it There are positive real numbers $c=c(g), \delta=\delta(g)$ 
such that if $x\in \mathcal{A}_g$ is a special point,
and $z\in F_g$ with $u(z)=x$ then
$$
H(z) \le c\Delta(x)^{\delta}.
$$
}

Daw-Orr \cite{DawOrr16} establish the generalization
for an arbitrary Shimura variety. It is
a result in "quantitative reduction theory".
%and a bit more
%technical to state as general Shimura varieties do not
%parameterize abelian varieties.

A few words about effectivity.
The point-counting approach to the Andr\'e-Oort 
conjecture had two sources of ineffectivity. 
First: the point-counting result in \cite{PilaWilkie} 
is in general ineffective, as mentioned above. 
%Indeed, in the 
%generality of an arbitrary o-minimal structure it cannot 
%reasonably be made effective \cite{PilaQJM}. 
But it has been 
made effective by Binyamini \cite{BinyaminiLogNoetherian}
in the structure 
$\mathbb{R}_{{\rm LN}\ \exp}$ of ``Log-Noetherian
functions'' (with exponentiation)
introduced there, which accommodates all 
the sets required in applying point-counting 
to Andr\'e-Oort.
The second source of ineffectivity was in the 
lower-bound for the size of Galois orbits
of special points, which depended on 
the famously ineffective Landau-Siegel lower bound for the 
class numbers of imaginary quadratic
fields. 
Until recently, the only cases of Andr\'e-Oort
known effectively were $\mathcal{A}_1^2$, due
independently to K\"uhne \cite{KuhneEffective, KuhneEffective2}
and Bilu-Masser-Zannier \cite{BiluMZ13},
and for non-compact curves in Hilbert modular surfaces
(and somewhat more generally) by Binyamini-Masser
\cite{BinyaminiMasser21}.

But the lower bounds have recently been effectivized 
in general by OpenAI \cite{OpenAIQRH}, in which the 
effective non-existence
of Landau-Siegel zeros is proved, see e.g. 
the human exposition \cite{Green}.
The earlier Edixhoven-Klingler-Ullmo-Yafaev
approach to Andr\'e-Oort, conditional on GRH,
also required good lower bounds for Galois orbits
of special points and, in addition, suitable upper 
bounds for small split primes in CM fields. 
These have now also been shown to follow from
QRH, see \cite{FreiSofosWidmer}.

\section{Schanuel's conjecture and Ax-Schanuel.}

We digress to introduce the classical Schanuel
conjecture on the transcendence properties of the
complex exponential, and the functional version proved
by Ax \cite{Ax71}.

Let $K\subset L$ be fields.
By ${\rm tr. deg.}_KL$ we denote the transcendence
degree of $L$ over $K$, i.e. the cardinality of
a transcendence basis.

\medbreak
\noindent
{\bf 7.1. Conjecture.\/} (Schanuel; 
see \cite{LangITN}, p30-31)
Let $z_1,\ldots, z_n\in\mathbb{C}$ be linearly
independent over $\mathbb{Q}$. Then
$${\rm tr. deg.}_{\mathbb{Q}}
\mathbb{Q}(z_1,\ldots, z_n, e^{z_1},\ldots, e^{z_n})
\ge n.$$

\medbreak

Schanuel's formulation beautifully encapsulates 
%arguably all (see \cite{???}) 
the arithmetic properties of the complex exponential. 
It is open in general for $n\ge 2$.
The special case where all the $z_i\in\overline{\mathbb{Q}}$ 
recovers the Lindemann-Weierstrass Theorem
(sometimes referred to as Lindemann's Theorem).

\medbreak
\noindent
{\bf 7.2. Theorem.} (Lindemann-Weierstrass, 
see e.g. \cite{LangITN})
{\it Let $z_1,\ldots, z_n\in\overline{\mathbb{Q}}$ 
be linearly independent over $\mathbb{Q}$. Then
$e^{z_1},\ldots, e^{z_n}$ are algebraically
independent over $\mathbb{Q}$.\/}

\medbreak

Functional versions, also (according to \cite{Ax71}) 
formulated by Schanuel, were proved by 
Ax \cite{Ax71}. A number of variants are given 
there, including a version for formal power series,
and a version in differential field. It is convenient
here to state an equivalent version in the complex
setting.
Let $u: \mathbb{C}^n \rightarrow (\mathbb{C}^\times)^n$
be the $n$th cartesian product of the
complex exponential, with
$D_n\subset \mathbb{C}^n \times (\mathbb{C}^\times)^n$
its graph.

\medbreak
\noindent
{\bf 7.3. Theorem.\/} 
(Classical Ax-Schanuel; 
complex version \cite{Tsimerman15})
{\it Let 
$V \subset \mathbb{C}^n \times (\mathbb{C}^\times)^n$ 
be an irreducible algebraic subvariety, 
and let $U$ be a connected, irreducible component 
of $V\cap D_n$. Then
$$
\dim_{\mathbb{C}} V \ge \dim_{\mathbb{C}}U +n
$$
unless the projection of $U$ to $\mathbb{C}^n$
is contained in a translate of a strict
$\mathbb{Q}$-linear subspace.}

\medbreak

If the inequality fails, it means that the dimension of the 
intersection component $U$ is atypically large;
its co-dimension in 
$\mathbb{C}^n \times (\mathbb{C}^\times)^n$
``should be'' ${\rm codim} D_n+{\rm codim V}$. The
theorem asserts that such atypical intersections 
are always ``explained'' by the intervention of
bi-algebraic subvarieties
(the translate of a $\mathbb{Q}$-linear
subspace).
The corresponding result for semi-abelian varieties is due
to Ax \cite{Ax72}. % (see also \cite{Kirby}), 
His argument in fact applies more generally to 
``almost semi-abelian varieties''
\cite{Bertrand} and is
extended to all commutative algebraic groups
in \cite{PilaTsimerman22} (with an application to the
rank one Riemann-Hilbert correspondence 
in \cite{PilaTsimerman26}).

\medbreak

The following special case, addressing the image
of an algebraic variety under exponentiation,
is the classical ``Ax-Lindemann''.

\medbreak
\noindent
{\bf 7.4. Theorem.\/}
(Classical Ax-Lindemann)
{\it Let $W\subset \mathbb{C}^n$ be a complex
algebraic variety. Then $u(W)$ is Zariski-dense
in $(\mathbb{C}^\times)^n$ unless $W$ is
contained in the translate of a strict 
$\mathbb{Q}$-linear subspace. Consequently, if
$u(W)\subset V$ for  an algebraic variety 
$V\subset (\mathbb{C}^\times)^n$ then there is
a translate $L$ of a $\mathbb{Q}$-linear subspace
such that $W\subset L$ and $u(L)\subset V$.}

\section{Functional transcendence: Ax-Lindemann.}

The functional transcendence question that arises in
applying the point-counting strategy for
Andr\'e-Oort to a Shimura 
variety $S$ is the following:

\smallskip

{\it For which algebraic subvarieties $W\subset \mathbb{H}_S$
is $u(W)$ not Zariski-dense in $S$; specifically,
when is $u(W)\subset V$?}

\smallskip

The answer is a direct generalization of the classical
Ax-Lindemann above; such containment is always
explained by the intervention of a weakly special 
(``bi-algebraic'') subvariety.
The results were proved in a sequence of papers,
beginning with the case of cartesian products of
the modular function \cite{PilaAO},
continuing with \cite{PilaTsimerman13} for
$\mathcal{A}_2$, and the compact case due to
Ullmo-Yafaev \cite{UllmoYafaevHyperbolic}.
The case of $\mathcal{A}_g$, crucial for the 
confirmation of Andr\'e-Oort for it later by 
Jacob, gives the form of all these statements.

\medskip
\noindent
{\bf 8.1. Theorem.\/} 
(Ax-Lindemann for $\mathcal{A}_g$; Pila-Tsimerman
\cite{PilaTsimerman14})
{\it Let $S={A}_g$ be the Siegel modular variety, 
let $\mathbb{H}_S=\mathbb{H}_g$ be the Siegel upper 
half-space, and let $u: \mathbb{H}_S\rightarrow S$
be the uniformizing map.
Let $W\subset \mathbb{H}_S, V\subset S$
be algebraic varieties. 

If $u(W)\subset V$ then there
is a weakly special subvariety $L\subset\mathbb{H}_S$
such that $W\subset L$ and $u(L)\subset V$.}

\medbreak

The general case for Shimura varieties 
(and even arithmetic varieties) is due to
Klingler-Ullmo-Yafaev 
\cite{KlinglerUllmoYafaevHyperbolicAL} who,
in the same paper, established the requisite
definability property. The extension to
mixed Shimura varieties is due to Gao 
\cite{GaoTowardsMixed}.
%{GaoAboutMixed,
The proofs of all these Ax-Lindemann results employ
point-counting. 

\section{Functional transcendence: 
Ax-Schanuel generalizations.}

Though not required for Andr\'e-Oort, it is natural to seek 
to generalize the full Ax-Schanuel statement.
Beyond their being natural statements, they have found
several applications. They are required in point-counting
approaches to more general cases of the Zilber-Pink
conjecture, in which they play a similar role
to Ax-Lindemann in Andr\'e-Oort. They have found further
application in sparsity results for integer points on
moduli spaces \cite{LawrenceVenkatesh, LawrenceSawin,
BaldiKlinglerUllmo}, in the Chabauty-Kim method 
\cite{Hast}, and in studying the degeneracy of the 
Betti map \cite{AndreCorvajaZannier, GaoBetti} 
(which plays a crucial role in uniformity
results related to Mordell-Lang 
\cite{GaoHabeggerHeights}).
They are employed in \cite{GMST24} to identify the
largest possible dimension of a
compact subvariety of $\mathcal{A}_g$, and in
\cite{PilaTsimerman22} to study problems generalizing 
elementary integrability.

As with Ax-Lindemann, Ax-Schanuel results have been 
established in a sequence of works by various authors 
in increasing generality.
The general result for Shimura varieties is due
to Mok-Pila-Tsimerman \cite{MokPilaTsimerman19}.
This was generalized to the mixed case 
of the universal abelian variety is by Gao 
\cite{GaoMixed} (see also \cite{GaoEnlarged}).

\medbreak
\noindent
{\bf 9.1. Theorem.\/} (Ax-Schanuel for Shimura varieties;
\cite{MokPilaTsimerman19})
{\it Let $u: \mathbb{H}_S\rightarrow S$ be the uniformizing
map of a Shimura variety, and $D\subset \mathbb{H}_S\times S$
its graph. Let $V \subset \mathbb{H}_S\times S$ 
be an irreducible algebraic subvariety, 
and let $U$ be a connected, irreducible component 
of $V\cap D_S$. Then
$$
\dim_{\mathbb{C}} V \ge \dim_{\mathbb{C}}U +\dim S
$$
unless the projection of $U$ to $S$
is contained in a strict
weakly special subvariety.}

\medbreak

A generalization to variations of pure polarized integral
Hodge structure is due to Bakker-Tsimerman
\cite{BakkerTsimerman} (see also
\cite{BakkerTsimermanLectures}).
%and to variations
%of mixed Hodge structure by Gao-Klingler
%\cite{GaoKlingler} and Chiu \cite{Chiu}.
In the Hodge-theoretic setting, one does not have
a uniformization map ``down'' to a target algebraic 
variety, but rather an period map $\phi$ ``up'' from 
the base $X$ of the variation into a period space
(see \S13 below).
%\Gamma\backslash D$, where the period domain $D$ is
%open in an algebraic variety $\hat{D}$
The natural ambient here is the quotient $\Gamma\backslash D$
of the associated ``weak Mumford-Tate domain'' $D$
by the image $\Gamma$ of the monodromy representation
of $\pi_1(X)$,
the conclusion is that atypical intersection
components are accounted for 
in a natural way, here 
by proper ``weak Mumford-Tate subvarieties''.

\medbreak
\noindent
{\bf 9.2. Theorem.\/} (Ax–Schanuel for variations of 
Hodge structures; \cite{BakkerTsimerman}). 
{\it In the above setup, let
$V \subset X \times \widehat{D}$ be an algebraic 
subvariety, let $X_D\subset X\times D$ be the fibre
product obtained from the quotient map
$D\rightarrow\Gamma\backslash D$ and the period
map $\phi: X\rightarrow \Gamma\backslash D$.
Let U be an irreducible analytic component of 
$V \cap X_D$ such that
$$
{\rm codim}_{X\times D} U < 
{\rm codim}_{X\times \widehat{D}} V  + 
{\rm codim}_{X\times D} X_D.
$$
Then the projection of U to X is contained in a 
strict weak Mumford–Tate subvariety.}

\medbreak

The method follows that of \cite{MokPilaTsimerman19},
but there are serious difficulties in finding a
suitable fundamental set, establishing definability
properties, and obtaining volume estimates.
The generalization to variations of mixed Hodge structures is
due to Gao-Klingler \cite{GaoKlingler} and to
Chiu \cite{Chiu} (and with derivatives
\cite{ChiuDerivatives}).

At this juncture, the approach to proving such
theorems via point-counting was supplanted by a 
more natural (and general) technology due to 
Blazquez-Sanz--Casale--Freitag--Nagloo
\cite{BSCFN21}. 
A very general Ax-Schanuel result along the lines of those 
above is a functional version of the Andr\'e-Grothendieck
period conjecture (with proofs via both techniques),
due to Bakker-Tsimerman \cite{BakkerTsimerman25}.

\section{Galois orbits of special points.}

A Galois lower-bound in a form suitable to prove the 
Andr\'e-Oort conjecture for $\mathcal{A}_1^n$ via 
point-counting is provided by one of the most classical 
results in number theory.

If an elliptic curve has complex multiplication then
its endomorphism ring is an order in a quadratic
imaginary number-field. (In the non-CM case,
the endomorphism ring is $\mathbb{Z}$.)
Let $\sigma$ be a singular modulus associated to
the order $\mathcal{O}$ of discriminant $\Delta$
(the discriminant uniquely characterizes the order).
It is classical that, in this case, all the CM points
for $\mathcal{O}$ are Galois conjugates, and that the 
Galois orbit is in bijection with the class group
of the order. In case of the maximal order in a given
field, i.e. its ring of integers, this is the class
number of the field. By Landau-Siegel 
(and now effective)
%by the famously ineffective
%(but now effective) Landau-Siegel bound 
%for class numbers
$$
[\mathbb{Q}(\sigma):\mathbb{Q}] 
\ge c(\epsilon) \Delta(\sigma)^{1/2-\epsilon},
$$
where $\Delta(\sigma)$ is the absolute value 
of the discriminant of $\mathcal{O}$.

But for other Shimura varieties the connection with
class numbers is more complicated, and Galois orbit bounds
harder to extract. Already for $\mathcal{A}_g$ for
higher $g$, the Galoius action on CM points is via
the reflex type norm, so the required estimate is about the
size of the image of a map between class groups.
A bound of the required form was raised earlier
as a question by Edixhoven \cite{EdixhovenMoonenOort}.

The required bounds were proved by Ullmo-Yafaev 
\cite{UllmoYafaevNombre} for $g\le 3$
and Tsimerman \cite{Tsimerman12} up to $g=6$
(and in both cases for all $g$ under GRH).
This result completed the proof of Andr\'e-Oort
for $\mathcal{A}_g^n$ for $g\le 6$
\cite{UllmoQuelques}.
But the argument requires
a combinatorial analysis whose complexity increases
with $g$.

\section{Heights and Galois orbits 
of special points.}

Going further, height bounds for special points
play an essential role. This reprises the role of height
bounds in the counterpart problems in group variety 
settings. 

Thus, in their work on unlikely
intersections for curves in $\mathbb{G}_{\rm m}^n$,
Bombieri-Masser-Zannier \cite{BMZ99}
crucially proved a bounded height
statement for the ``just likely'' intersections.
In their subsequent work \cite{BMZ07}, 
they framed a general ``Bounded Height Conjecture''
in connection with their version of the
Zilber-Pink conjecture.
The Bounded Height Conjecture was proved by Habegger,
both in the multiplicative setting \cite{HabeggerBHC}
and that of abelian varieties \cite{HabeggerAbelian}.
These results (until OpenAI \cite{OpenAIAbelianZP}
claimed a full proof of  the Zilber-Pink conjecture
for abelian varieties over $\overline{\mathbb{Q}}$) 
afforded the most general result
on the Zilber-Pink conjecture in these settings.
The failure of ``bounded height'' in the setting of
Shimura varieties was noted in \cite{BMZ07}
(indeed, already singular moduli, unlike torsion points
of semi-abelian varieties, 
do not have bounded height).
A weakly bounded height conjecture for 
``likely intersections'' is framed
in \cite{HabeggerWeakly}.

The key result required for Andr\'e-Oort is that the 
{\bf logarithmic\/} Weil height (i.e. $h(x)=\log H(x)$
where $H(x)$ is the multiplicative 
Weil height introduced above) of a {\bf special point\/}
$x\in \mathcal{A}_g$ is ``discriminant negligible'', 
i.e. that
$$
h(x) \le c(g, \epsilon) \Delta(x)^\epsilon
$$
for some $c(g,\epsilon)$, for all positive $\epsilon$.

Tsimerman \cite{Tsimerman18} observes that the above
follows from the just-proved Averaged Colmez Conjecture.
Colmez \cite{Colmez93} formulates a precise conjecture for
the Faltings height of a CM abelian variety. 
This remains open in general but its averaged
version was proved independently by
Andreatta-Goren-Howard-Madapusi Pera
\cite{AGHMP18} and
Yuan-Zhang \cite{YuanZhang}. (The average is discriminant
negligible, while the Faltings height is bounded below
over $\mathcal{A}_g$; \cite{Bost}, for an explicit 
bound see \cite{Pazuki}).
A clever argument using the isogeny estimates of
Masser-W\"ustholz then shows that the Galois orbit of
a special point must be at least of size a positive
power of its discriminant.

\medbreak
\noindent
{\bf 11.1. Theorem.} (\cite{Tsimerman18})
{\it Let $g$ be a positive integer.
There exists $c_g, \delta_g > 0$ such if $\Phi$ is a
primitive CM type for a CM field $E$, and if $A$ 
is an abelian variety of dimension $g$ with endomorphism
ring equal to the full ring of integers $\mathcal{O}_E$ 
and CM type $\Phi$, then the
field of moduli $\mathbb{Q}(A)$ of $A$ satisfies
$$
[\mathbb{Q}(A) : Q] \ge c_g |{\rm Disc}(E)|^{\delta_g}.
$$
}

\medbreak

From this one deduces a bound for a general CM point, and 
this gives the final ingredient to establish the
Andr\'e-Oort conjecture for $\mathcal{A}_g$. 

\medbreak
\noindent
{\bf 11.2. Theorem.\/} (Tsimerman \cite{Tsimerman18})
{\it The Andr\'e-Oort conjecture holds
for $\mathcal{A}_g$.}

\medbreak

For general
Shimura varieties, we have seen that the required ingredients
1, 2, 3 have been obtained. Only the Galois orbit lower
bound for special points on general Shimura varieties
was outstanding as of 2018.
Though formally not part of the work cited by the Fields Medal
Committee, we mention the further developments leading
to the announcement \cite{PilaShankarTsimerman21}
of a full proof of the Andr\'e-Oort conjecture.

For a general Shimura variety, its points do not
parameterize abelian varieties (or even, to present knowledge,
algebraic varieties). Two tools crucially employed in 
obtaining the Galois orbit lower bound are therefore absent.
Namely, one lacks the Faltings height, and 
isogeny estimates.

The absence of isogeny estimates was circumvented by
Binyamini-Schmidt-Yafaev. Special points give rise to
algebraic points on the graph $D_S$ of the uniformization
$u: \mathbb{H}_S\rightarrow S$.
They no longer have bounded degree. By an idea of Schmidt,
implemented in \cite{Schmidt} for torsion points, 
a sufficiently sharp (poly-logarithmic) counting
bound for algebraic points in height and degree
aspect on $D_S$  (restricted to
a fundamental domain, of course), combined with the fact that
there are known to be ``many'' such points, can be used to force
a degree lower bound. (The algebraic part of $D_S$ is empty,
so the degrees and heights must 
respect the upper bound.)
A sharp counting result applicable in the 
general Shimura variety context
is due to Binyamini \cite{BinyaminiLogNoetherian}, and
with Yafaev \cite{BSY23} this is
implemented to show that if the heights of special points
on $S$ are ``discriminant negligible'' then a suitable
Galois orbit lower bound follows.

The paper \cite{PilaShankarTsimerman21} then establishes
the height bound. This uses the fact that certain special points
on a general Shimura variety can be related to special points
on $\mathcal{A}_g$, and by careful control one can use this
to construct an ``almost canonical'' height on special points
on an arbitrary Shimura variety which is 
discriminant negligible. The height bound again ultimately 
rests on the proof of the Averaged Colmez 
conjecture.

\section{O-minimal complex geometry.}

In a series of papers, Peterzil and Starchenko develop
complex analysis in the setting of an o-minimal structure 
\cite{PeterzilStarchenkoExpansions, 
PeterzilStarchenkoExpansions2,
PeterzilStarchenkoNonStandard,
PeterzilStarchenkoCategories}
and the survey \cite{PeterzilStarchenkoTame}. 

It is a beautiful and remarkable development in which
winding numbers replace the use of power series,
which are unavailable in general:
the universe of an o-minimal structure is not, 
in general, the real numbers, but a real closed field
which may be non-archimedean (and not locally compact).
In the one-variable setting one nevertheless recover
all the standard properties (maximum principle, open mapping,
Cauchy-Riemann relations and infinite differentiability,
identity principle).

In several variables, Peterzil and Starchenko
establish o-minimal analogues of several
key classical results, including
Chow's Theorem, Remmert-Stein Theorem, and
Removable Singularities. The o-minimal analogues 
in general include the classical statements in case 
the ambient structure is $\mathbb{R}_{\rm an}$
and, in effect, replace properness assumptions by
definability.

\medbreak
\noindent
{\bf 12.1. Theorem.\/} (``Definable Remmert-Stein''; 
\cite{PeterzilStarchenkoCategories}. See also
\cite{PeterzilStarchenkoNonStandard}). 
{\it Let $M$ be a definable complex manifold and $E$ 
a definable complex analytic subset of $M$.
If $A$ is a definable, complex analytic subset of $M-E$, 
then its closure $\overline{A}$ is a
complex analytic subset of $M$.\/}

\medbreak

The following is a slight generalization of the result 
stated by Peterzil-Starchenko.

\medbreak
\noindent
{\bf 12.2. Theorem.\/} (``Definable Chow''; 
\cite{PeterzilStarchenkoCategories},
\cite{MokPilaTsimerman19})
{\it Let $Y$ be a quasi-projective algebraic
variety, and let $A \subset Y$ be definable, 
complex analytic, and closed in $Y$. Then
$A$ is algebraic.\/}

\medbreak

These powerful statements are applicable, by the 
definability results above, and those further below,
to the (suitably restricted)
uniformization maps of Shimura varieties, and more 
generally to the period maps of variations of Hodge
structure.

\section{Hodge theory.}

The primary object of Hodge theory is the 
{\it Hodge structure\/} associated to a suitable complex 
manifold. It records the linear algebraic data
associated to the interplay of topological data
and the complex structure.
We will consider only the case of smooth
and projective complex
algebraic varieties.

The Betti cohomology groups, tensored up to 
$\mathbb{C}$, retain an integral structure. The
smooth De Rham cohomology groups $H^k(X,\mathbb{C})$
of such a variety $X$ have 
a {\it Hodge decomposition\/}, a direct sum
decomposition into subspaces $H^{p,q}$ over
non-negative integers $p,q$ with $p+q=k$,
spanned by forms of type $(p,q)$, meaning locally 
represented as sums over wedge products of 
$p$ of the $dz_i$ and $q$ of the $d\overline{z}_j$.
%in degree $p$, anti-holomorphic in degree $q$.
These satisfy that
$H^{p,q}$ is the complex conjugate of $H^{q,p}$.
This data gives a {\it weight $k$ $\mathbb{Z}$-Hodge
structure\/} $V_{\mathbb{Z}}$
associated to $X$. But there is further
a {\it polarization\/} $Q$ coming from the intersection
product (and the projective embedding), it
is a bilinear form on $H^k(X,\mathbb{C})$
enjoying suitable (anti-/)symmetry
and positivity properties.

The data associated with such a weight $k$ polarized 
$\mathbb{Z}$-Hodge structure,  $(V_{\mathbb{Z}}, Q)$,
is naturally parameterized by the 
associated {\it period domain\/}, an open subset $D$
of a suitable flag variety $\hat{D}$. While
$\hat{D}$ is complex projective variety, 
the subset $D$ is defined by positivity conditions.
The period domain is a quotient $G/M$
of $G={\bf G}(\mathbb{R})^+$, the connected
component of the identity of ${\bf G}(\mathbb{R})$
for some connected semi-simple linear algebraic group
$\mathbb{Q}$-group ${\bf G}$ by a compact subgroup $M$
(in general not maximal). The isomorphism types
of such Hodge structures
are realized by the {\it period space\/}, a quotient
$\Gamma\backslash D$ by some arithmetic group
$\Gamma$.

The prototype examples are the Siegel upper half-space. 
$$
\mathcal{A}_g={\rm Sp}_{2g}(\mathbb{Z})\backslash \mathbb{H}_g=
{\rm Sp}_{2g}(\mathbb{Z})\backslash {\rm Sp}_{2g}(\mathbb{R})/
U(g).
$$
But here the
compact subgroup $U(g)$ is maximal and the quotient 
${\rm Sp}_{2g}(\mathbb{Z})\backslash \mathbb{H}_g$
is an algebraic variety, 
in fact a Shimura variety (which all arise
in this way); but in general the period space
has no algebraic structure \cite{GriffithsRoblesToledo}.
Every smooth projective complex
algebraic variety $X$ of the same topological type
will give a polarized weight $k$ $\mathbb{Z}$-Hodge structure,
for $k$ in the appropriate range,
corresponding to a point in this space. 

The basic object for our discussion is the 
``variation of Hodge structure'', developed by Griffiths 
\cite{GriffithsPeriodsI, GriffithsPeriodsII, 
GriffithsPeriodsIII} and the survey
\cite{Griffiths70}. 
Suppose that $\mathcal{X}\rightarrow S$ is a morphism
of algebraic varieties whose generic fibre is a smooth
and projective. The base $S$ is assumed to be
quasi-projective; by shrinking it we may assume that
all fibres are smooth and projective.
The fibres $\mathcal{X}_s, s\in S$ are then homeomorphic 
as real manifolds, and we may think of
$\mathcal{X}\rightarrow S$ 
as a varying family of complex structures 
on an underlying smooth compact real manifold, giving
a family $V\rightarrow S$ of polarized 
$\mathbb{Z}$-Hodge structures. These give a bundle over 
the base with a holomorphically varying filtration, 
with regular, flat connection, and satisfying 
Griffiths transversality. They may then be considered
abstractly.
Variations arising from families $\mathcal{X}\rightarrow S$
are called {\it geometric\/}.

For example, if we take the universal family of abelian
varieties $\mathcal{X}_g\rightarrow\mathcal{A}_g$
over $\mathcal{A}_g$ (with some suitable level
structure that we will suppress), then the
periods, as multivalued functions on the base, map
$\mathcal{A}_g$ to the Siegel upper-half space
$\mathbb{H}_g$. The monodromy is (some congruence
subgroup, depending on the level structure, of)
${\rm Sp}_{2g}(\mathbb{Z})$, and passing to the quotient
gives a prototypical variation of Hodge structure.
Though prototypical, this variation is 
again extremely special,
as the period map is onto, the quotient map
is a biholomorphism of the quotient domain
with $\mathcal{A}_g$. (With the level structure:
on $\mathcal{A}_g$ itself the map from $\mathbb{H}_g$
is ramified.) This property is enjoyed by all 
Shimura varieties (though in the absence of the 
Hodge conjecture they are not all 
known to parameterize families of algebraic varieties).

In general, the family $\mathcal{X}\rightarrow S$ 
gives rise to a period map into the period domain. 
The map is in general multi-valued if $S$ is not
simply connected; quotienting
by the monodromy one gets the {\it period map\/} into the 
period space $\Gamma\backslash D$. 
In general, the period map
is not onto, but of high codimension.
And, in general, the period space has no algebraic 
structure.
Nevertheless, Griffiths conjectured (in the proper case)
\cite{Griffiths70} that the image 
of the period map is algebraic, and further
properties, including that the 
Griffiths bundle is semi-ample.

\medbreak
\noindent
{\bf 13.1. Conjecture.} 
(Griffiths, \cite{Griffiths70} (10.5)
``Existence of automorphic forms'')
Given $ \phi: S \rightarrow \Gamma\backslash D$ as above, 
the image $\phi(S)$ of $S$ is an algebraic variety.

\medbreak

Griffiths shows the conjecture holds when the image is compact,
and Sommese later affirmed it when the singularities are
isolated.

The {\it Hodge locus\/}
${\rm HL}(S, V) \subset S$ associated to a 
polarized variation of Hodge structures
$V$ is the set of points $s \in S$ for which exceptional 
Hodge tensors (i.e. $\mathbb{Z}$-classes with pure type 
with respect to the Hodge decomposition) occur for $V_s$. 
The locus ${\rm HL}(S, V)$
is easily seen to be a countable union of irreducible 
complex analytic subvarieties of $S$, called {\it special 
subvarieties\/} of $S$ associated to $V$. 
If $V = R^kf_{*}\mathbb{Q}$ for 
a smooth proper morphism  $f : \mathcal{X} \rightarrow S$,
the Hodge conjecture implies that
the exceptional Hodge tensors in $V_s$ come from 
exceptional algebraic cycles in some product 
$\mathcal{X}_s^N$.
One can then deduce that every special subvariety
(component of the Hodge locus)
of $S$ is algebraic. 
A seminal result was the affirmation
of this property unconditionally, by
Cattani-Deligne-Kaplan \cite{CattaniDeligneKaplan},
which is viewed as significant evidence for the
Hodge conjecture.

\medbreak
\noindent
{\bf 13.2. Theorem.} 
(Cattani-Deligne-Kaplan \cite{CattaniDeligneKaplan})
{\it The special subvarieties of $S$ associated to $V$ 
are algebraic, i.e. the Hodge locus ${\rm HL}(S, V)$ 
is a countable union of closed irreducible algebraic 
subvarieties of $S$.\/}

\medbreak

If $S$ is a Shimura variety, the Hodge locus
is precisely the locus of special subvarieties, and is
(generally) a countably infinite
collection of algebraic subvarieties defined 
over $\overline{\mathbb{Q}}$.
Whether the components of the Hodge locus are in general
defined over $\overline{\mathbb{Q}}$ 
(when the variation is so defined) remains open.
By \cite{KlinglerOtwinowskaUrbanik} 
it would follow in general 
from the result for zero-dimensional components. 
Thus while the framework of
``unlikely intersections'' and the Zilber-Pink
conjecture (including the Andr\'e-Oort conjecture)
have been generalized to the Hodge loci of
(graded polarizable, admissible) variations of 
mixed Hodge structure over a smooth complex
quasi-projective base
\cite{KlinglerHodgeLoci},
results to date address mainly geometric
aspects, in which Ax-Schanuel results play a 
prominent role.

\section{O-minimality and Hodge theory.}

The potential utility of o-minimal methods in
Hodge theory was highlighted very early in its 
development, as a setting that can accommodate maps 
such as the logarithm that are neither algebraic nor 
subanalytic at the singularity: see the remarks in
\cite{DriesMillerCategories} and \cite{SchmidVilonen}.
However, it seems that deeper applications 
awaited the development of complex geometry in these
settings by Peterzil-Starchenko, their
``definable Chow,'' and its profound elaboration 
into ``o-minimal GAGA'', by Jacob and his co-workers.

First we return to definability issues.
The period space is a definable 
(in $\mathbb{R}_{\rm alg}$) manifold,
but the map into it given by a variation of
Hodge structure is highly transcendental.
In \cite{BakkerKlinglerTsimerman}, 
Bakker-Klingler-Tsimerman show that period map is definable
(in $\mathbb{R}_{\rm an\ exp}$).

\medbreak
\noindent
{\bf 14.1. Theorem.\/} 
(Bakker-Klingler-Tsimerman \cite{BakkerKlinglerTsimerman})
{\it Let $V$ be a polarized variation of pure Hodge structures
of weight $k$ over
a smooth complex quasi-projective variety $S$. 
Let  $\Phi_S : S \rightarrow \Gamma\backslash D$
be the period map associated to $V$,
where $\Gamma$ is a (neat) arithmetic group 
containing the monodromy.
Then $\Phi_S$ is 
definable in $\mathbb{R}_{\rm an\ exp}$.\/}

\medbreak

It is a delicate result that
depends crucially on the Schmid ${\rm SL}_2$-orbit
theorem \cite{SchmidVariation} and
asymptotics for Hodge norms.
As an immediate corollary of the definability result
and Peterzil-Starchenko’s ``definable Chow'' Theorem,
they obtain an alternative proof
of the algebraicity of the Hodge locus
(Cattani-Deligne-Kaplan theorem).
The definability and quasi-projectivity of Hodge 
loci components
%\cite{BakkerBrunebarbeTsimerman23}
are extended to mixed period mappings 
by Bakker-Brunebarbe-Klingler-Tsimerman
\cite{BakkerBrunebarbeKlinglerTsimerman24}
(recovering the extension of Cattani-Deligne-Kaplan to 
this setting due to Brosnan-Pearlstein-Schnell 
\cite{BrosnanPearlsteinSchnell}).

The proof of the Griffiths conjecture required the 
development of a GAGA principle in the o-minimal 
context, an analogue of the seminal GAGA of Serre
\cite{SerreGAGA}. This is carried out by 
Bakker-Brunebarbe-Tsimerman
in \cite{BakkerBrunebarbeTsimerman22},
and provides a flexible and powerful tool
in the o-minimal context.

\medbreak
\noindent
{\bf 14.2. Theorem.} (Bakker-Brunebarbe-Tsimerman
\cite{BakkerBrunebarbeTsimerman22}). 
{\it Let $\mathcal{R}$ be an o-minimal 
structure over the real field.
Let $X$ be a separated algebraic space of finite type over
$\mathbb{C}$, and $X^{\rm def}$ the
associated definable analytic space. 
Then the {\rm definablization functor\/}
$(-)^{\rm def} :
{\rm Coh}(X) \rightarrow  {\rm Coh}(X^{\rm def})$
is fully faithful and exact, and its essential image is
closed under sub-objects and quotients.

Thus, in particular, definable coherent subsheaves of 
algebraic coherent
sheaves are algebraic.\/}

\medbreak

Note that $X$ is not required 
to be proper over $\mathbb{C}$ but that, in contrast
to Serre’s classical GAGA, it is not true that 
every definable coherent sheaf is algebraic.
%Note further that $(-)^{\rm def}$ is in general not 
%essentially surjective. 
The extension to the mixed case is in
\cite{BakkerBrunebarbeTsimerman23}.

Combining o-minimal GAGA with 
algebraization theorems of
Artin, the authors show that proper definable images of 
algebraic varieties are algebraic.
This enables them to deduce that period images are
algebraic, confirming the conjecture 
of Griffiths.

\medbreak
\noindent
{\bf 14.3. Theorem.\/} (Bakker-Brunebarbe-Tsimerman
\cite{BakkerBrunebarbeTsimerman22}) 
{\it Let $X$ be a reduced separated
algebraic space of finite type over $\mathbb{C}$
and $\phi : X^{\rm an} \rightarrow \Gamma\backslash D$ 
a period map. Then

(1) $\phi$ factors (uniquely) as 
$\phi = \iota \circ f^{\rm an}$ where $f : X \rightarrow Y$
is a dominant map
of (reduced) finite-type algebraic spaces and 
$\iota : Y^{\rm an} \rightarrow \Gamma\backslash D$ 
is a closed
immersion of analytic spaces;

(2) the Griffiths $\mathbb{Q}$-bundle $L$ restricted to $Y$ 
is the analytification of an ample
algebraic $\mathbb{Q}$-bundle, and in particular $Y$ is a 
quasi-projective variety.\/}

\medbreak

Beyond clarifying foundational aspects
of variational Hodge theory, the results have led
to further important progress.
These include the construction of
toroidal (Deng-Tsimerman \cite{DengTsimerman25}) and
Baily-Borel 
(Bakker-Filipazzi-Mauri-Tsimerman \cite{BFMT25})
compactifications 
of period images, 
with associated ampleness results,
both crucial objects in their classical incarnation
for arithmetic quotients, with the latter having
important applications to birational geometry via
the $b$-semiampleness conjecture,
and the affirmation
of the linear Shafarevich conjecture and the
construction of algebraic Shafarevich morphisms 
for arbitrary varieties (Bakker-Brunebarbe-Tsimerman
\cite{BakkerBrunebarbeTsimerman24a}).

\section{Other work.} I have not touched on many 
other beautiful results of Jacob that are not 
directly connected to the work in the citation.
See arXiv and MatSciNet.

\section{Conclusion.}

Jacob's work has indeed established o-minimality as 
an essential setting for Hodge theory and 
arithmetic geometry: Hodge theory is tame, and the
consequences are profound.

%He has further (with co-authors) 
%developed o-minimal GAGA as a foundational tool,
%and made decisive advances in our understanding of
%the basic objects of Hodge theory. He has, 
%along the way, made spectacular advances in 
%our understanding of the arithmetic of special 
%points of Shimura varieties.

%\medbreak
%\noindent
%{\bf AI disclosure.\/} AI was only used as it
%appeared in Google searches on bibliographical and 
%terminological issues. All em-dashes, semi-colons
%and errors are the sole responsibility of the author.

\medbreak
\noindent
{\bf Acknowledgements.\/} 
Many thanks to Ben Bakker and Jacob Tsimerman for their 
comments on and corrections to an earlier draft. 
All remaining errors and infelicities are
the responsibility of the author. This paper was written
while I was a visitor, first at the Institute
for Advanced Study in Princeton and then at the 
Weizmann Institute of Science in Rehovot. 
I am very grateful to both institutes for their
warm hospitality, superb environment, and generous
financial support. At the IAS I was supported
by the  Charles Simonyi Endowment, and at the WIS
by a research grant from Magnus Konow in honour 
of his mother Olga Konow Rappaport, and by the
European Union (ERC, SharpOS, 101087910). 

\bigbreak

\end{document}